\documentclass{amsart}

\usepackage[normalem]{ulem}

\usepackage{lineno}

\usepackage{soul}

\usepackage{amssymb, amsmath, amsthm, hyperref, euscript, color}
\usepackage[dvipsnames]{xcolor}

\usepackage{amscd}

\usepackage{bbm}

\usepackage{enumitem}

\usepackage[english]{babel}

\numberwithin{equation}{section}

\theoremstyle{plain}

\newtheorem{theorem}{Theorem}

\newtheorem{remark}{Remark}

\newtheorem{lemma}{Lemma}

\newtheorem{thm}{Theorem}

\newtheorem{cor}[thm]{Corollary}

\newtheorem{prop}[thm]{Proposition}

\theoremstyle{definition}

\newtheorem{definition}{Definition}

\newcommand{\R}{\mathbb{R}}

\newcommand{\Z}{\mathbb{Z}}

\newcommand{\SO}{\mathrm{SO}}
\newcommand{\SU}{\mathrm{SU}}

\DeclareMathOperator{\Min}{min}

\DeclareMathOperator{\Sec}{sec}
\DeclareMathOperator{\Ric}{Ric}
\DeclareMathOperator{\Gr}{Gr}

\DeclareMathOperator{\Id}{Id}

\DeclareMathOperator{\Dist}{dis}

\begin{document}

\author{Boris Stupovski and Rafael Torres}

\title[Simply connected 5-manifolds of positive biorthogonal curvature]{Existence of Riemannian metrics with positive biorthogonal curvature on simply connected 5-manifolds}

\address{Scuola Internazionale Superiori di Studi Avanzati (SISSA)\\ Via Bonomea 265\\34136\\Trieste\\Italy}

\email{bstupovs@sissa.it}

\email{rtorres@sissa.it}

\subjclass[2010]{Primary 53C20, 53C21; Secondary 53B21}

\maketitle

\emph{Abstract}: Using recent work of Bettiol, we show that a first-order conformal deformation of Wilking's metric of almost-positive sectional curvature on $S^2\times S^3$ yields a family of metrics with strictly positive average of sectional curvatures of any pair of 2-planes that are separated by a minimal distance in the 2-Grassmanian. A result of Smale's allows us to conclude that every closed simply connected 5-manifold with torsion-free homology and trivial second Stiefel-Whitney class admits a Riemannian metric with a strictly positive average of sectional curvatures of any pair of orthogonal 2-planes. 

\section{Introduction and main results}\label{Introduction}

Let $(M, g)$ be a compact Riemannian $n$-manifold and let $\Sec_g$ be the sectional curvature of the metric. We often abuse notation and denote the Riemannian metric by $(M, g)$ as well. For each 2-plane\begin{equation}\sigma\in \Gr_2(T_pM) = \{X\wedge Y\in {\Lambda}^2 T_pM : ||X\wedge Y||^2 = 1\},\end{equation} let $\sigma^\perp\subset T_pM$ be its orthogonal complement. That is, there is a $g$-orthogonal direct sum decomposition $\sigma \oplus \sigma^{\perp} = T_pM$ at a point $p\in M$.

\begin{definition}\label{Definition Biorthogonal Curvature} The biorthogonal curvature of a 2-plane $\sigma\in \Gr_2(T_pM)$  is\begin{equation}\label{Equation 1}\Sec_g^{\perp}(\sigma):= \underset{\begin{subarray}{c}
\sigma'\in \Gr_2(T_pM)\\ 
\sigma' \subset \sigma^{\perp}\end{subarray}}{\Min}\frac{1}{2}(\Sec_g(\sigma) + \Sec_g(\sigma'))\end{equation}(cf. \cite[Section 5.4]{[Be2]}). We say that $(M, g)$ has positive biorthogonal curvature $\Sec_g^\perp > 0$ if (\ref{Equation 1}) is positive for every $\sigma\in \Gr_2(T_pM)$ at every point $p\in M$ .
\end{definition}

A stronger curvature condition is the following. Choose a distance on the Grassmanian bundle $\Gr_2(TM)$ that induces the standard topology.

\begin{definition}\label{Definition Distance Curvature} The distance curvature of a 2-plane $\sigma \subset T_pM$ is\begin{equation}\label{Equation 2}\Sec_g^{\theta}(\sigma):= \underset{\begin{subarray}{c}
\sigma'\in \Gr_2(T_pM)\\ 
\Dist(\sigma, \sigma')\geq \theta\end{subarray}}{\Min}\frac{1}{2}(\Sec_g(\sigma) + \Sec_g(\sigma'))\end{equation} for  each $\theta > 0$  (cf. \cite[Section 5.2]{[Be2]}). We say that $(M, g^\theta)$ has positive distance curvature $\Sec_{g^\theta} > 0$, if for every $\theta > 0$ there is a Riemannian metric $(M, g^\theta)$ for which (\ref{Equation 2}) is positive for every $\sigma\in \Gr_2(T_pM)$ at every point $p\in M$.
\end{definition}

Bettiol \cite{[Be3]} classified up to homeomorphism closed simply connected 4-manifolds that admit a Riemannian metric of positive biorthogonal curvature by constructing metrics of positive distance curvature on $S^2\times S^2$   \cite[Theorem, Proposition 5.1]{[Be1]}, \cite[Theorem 6.1]{[Be2]} and showing that positive biorthogonal curvature is a property that is closed under connected sums \cite[Proposition 7.11]{[Be2]}, \cite[Proposition 3.1]{[Be3]}.

In this paper, we extend Bettiol's results to dimension five. More precisely, we build upon Bettiol's work and show that an application of a first-order conformal deformation to Wilking's metric $(S^2\times S^3, g_W)$ of almost-positive sectional curvature \cite{[Wilking]} yields the main result of this note.

\begin{thm}\label{Theorem A} For every $\theta > 0$, there is a Riemannian metric $(S^2\times S^3, g^{\theta})$ such that

(a)  $\Sec_{g^\theta}^\theta > 0$.

(b) there is a limit metric $g^0$ such that $g^{\theta}\rightarrow g^0$ in the $C^k$-topology as $\theta\rightarrow 0$ for $k\geq 0$.

(c) $g^\theta$ is arbitrarily close to Wilking's metric $g_W$ of almost-positive curvature in the $C^k$-topology for $k\geq 0$.

(d) $\Ric_{g^{\theta}} > 0$.

(e) There is a 2-plane $\sigma\in \Gr_2(T_pS^2\times S^3)$ with $\Sec_{g^{\theta}}(\sigma) < 0$.

In particular, there is a Riemannian metric of positive biorthogonal curvature on $S^2\times S^3$. 
\end{thm}

The next corollary is a consequence of coupling Theorem \ref{Theorem A} with a  classification result of Smale \cite{[S]}. 

\begin{cor}\label{Corollary B} Every closed simply connected 5-manifold with torsion-free homology and zero second Stiefel-Whitney class admits a Riemannian metric of positive biorthogonal curvature.
\end{cor}

The hypothesis imposed on the homology and the second Stiefel-Whitney class of the manifolds of Corollary \ref{Corollary B} are technical in nature; cf. Remark \ref{Remark 2}. Indeed, an examination of the canonical metric on the Wu manifold yields the following proposition. 

\begin{prop}\label{Proposition Wu manifold} The symmetric space metric $(\SU(3)/\SO(3), g)$ has positive biorthogonal curvature. 
\end{prop}

The Wu manifold has second homology group of order two and non-trivial second Stiefel-Whitney class.

\subsection{Acknowledgements:} We thank the referee and Renato Bettiol for useful input that allowed us to significantly improve the paper. R. T. thanks Nicola Gigli for very helpful discussions during the production of the paper.

\section{Constructions of Riemannian metrics of positive biorthogonal curvature}\label{Section ConstructionsRMetrics} 

\subsection{Wilking's metric of almost-positive curvature on $S^2\times S^3$}\label{Section Wilking's Metric} We follow the exposition in \cite[Section 5]{[Wilking]} to describe Wilking's construction a metric of almost positive curvature on the product of projective spaces $\R P^2\times \R P^3$ and its pullback to $S^2\times S^3$ under the covering map; see \cite[Section 5]{[Z1]} for a discussion relating these two constructions. 

The unit tangent sphere bundle of the 3-sphere\begin{equation}\label{eq:101}
    T_1(S^3)=S^2\times S^3\,,
\end{equation}embeds into $\mathbb{R}^4\times \mathbb{R}^4=\mathbb{H}\times \mathbb{H}$ as a pair of orthogonal unit quaternions\begin{equation}\label{eq:102}
    S^3\times S^2 = \{(p,v)\in\mathbb{H}\times\mathbb{H}:|p|=|v|=1, \langle p,v \rangle = 0\}\subset \mathbb{H}\times \mathbb{H}\,,
\end{equation}
where $\langle x,y \rangle = \mathrm{Re}(\bar{x}y)$, $|x|^2=\langle x,x \rangle$ and $\bar{x}$ denotes quarternion conjugation of $x$. The group $G=\mathrm{Sp(1)}\times \mathrm{Sp(1)} \simeq S^3\times S^3$ acts on $S^2\times S^3$ by 
\begin{equation}\label{eq:103}
    (q_1,q_2)\star(p,v)=(q_1 p \bar{q}_2,q_1 v \bar{q}_2)\,,
\end{equation}
for $q_1, q_2 \in \mathrm{Sp(1)}$ and $(p,v) \in S^2 \times S^3$. This action is effectively free and transitive. The isotropy group of the point $(1,i)\in S^2\times S^3$ is\begin{equation}H=\{(e^{i\phi},e^{i\phi})\in \mathrm{Sp(1)}\times \mathrm{Sp(1)}\}\subset G.\end{equation} Thus, $S^2\times S^3\simeq G/H$ is a homogeneous space.

In order to put a metric on $S^2\times S^3$, Wilking first defines a left invariant metric $g$ on $G=\mathrm{Sp(1)}\times \mathrm{Sp(1)}$ as follows. Let\begin{equation}g_0((X,Y),(X',Y'))=\langle X, Y \rangle + \langle X', Y' \rangle\end{equation} for $(X,Y),(X',Y')\in \mathfrak{sp}(1)\oplus \mathfrak{sp}(1)=\mathrm{Im}(\mathbb{H})\oplus\mathrm{Im}(\mathbb{H})$, denote a bi-invariant metric. In terms of $g_0$, the metric $g$ is 
\begin{equation}\label{eq:104}
    g((X,Y),(X',Y'))=g_0(\Phi(X,Y),(X',Y'))\,,
\end{equation}
where $\Phi$ is a $g_0$-symmetric, positive definite endomorphism of $\mathfrak{sp}(1)\oplus \mathfrak{sp}(1)$ given by
\begin{equation}\label{eq:105}
    \Phi = \Id - \frac{1}{2}P\,,
\end{equation}
and $P$ is the $g_0$-orthogonal projection onto the diagonal subalgebra\begin{equation}\Delta\mathfrak{sp}(1)\subset \mathfrak{sp}(1)\oplus\mathfrak{sp}(1);\end{equation}see \cite[p. 125]{[Wilking]}.

Wilking's doubling trick guarantees the existence of a diffeomorphism 
\begin{equation}\label{eq:106}
    G/H \simeq \Delta G\backslash G\times G/\{1_G\}\times H\,,
\end{equation}
where $ \Delta G\backslash$ denotes the quotient by the left diagonal action of $G$ on $G \times G$ and $H$ acts on the second factor from the right. Consider the product $(G\times G, g + g)$ (cf. \eqref{eq:104}) and the induced metric on $S^2\times S^3 \simeq \Delta G\backslash G\times G/\{1_G\}\times H$ that we denote by $g_W$. That is, Wilking's metric $(S^2\times S^3, g_W)$ is the metric that makes the quotient submersion
\begin{equation}\label{eq:107}
    (G\times G, g\oplus g) \rightarrow (\Delta G\backslash G\times G/\{1_G\}\times H, g_W)
\end{equation}
into a Riemannian submersion. Wilking has shown that $(S^2 \times S^3, g_W)$ has almost positive curvature, with flat 2-planes located on two hypersurfaces. These hypersurfaces are both diffeomorphic to $S^2\times S^2$, and they intersect along an $\R P^3$ \cite[Corollary 3, Proposition 6]{[Wilking]}. However, except for points that lie on four disjoint copies of $S^2$ inside these two hypersurfaces, there is a unique flat 2-plane. At each point in these four 2-spheres, there is a one parameter family of flat 2-planes and neither the distance curvature nor the biorthogonal curvature of the metric $g_W$ are strictly positive at any of these points.

\section{Proofs}\label{Section Proofs}

\subsection{Proof of Theorem \ref{Theorem A}} We follow Bettiol's construction of metrics of positive distance curvature  on $S^2 \times S^2$ \cite[Theorem]{[Be1]}, \cite[Theorem 6.1]{[Be2]}, and apply a first-order conformal deformation to Wilking's metric $(S^2\times S^3, g_W)$ that was described in Section \ref{Section Wilking's Metric}. This yields metrics of positive distance curvature as in Definition \ref{Definition Distance Curvature}, which converge to a metric $g^0$ as $\theta$ tends to zero in the $C^k$-topology. 
 
\begin{definition}\label{Definition confdef}
    Let $(M,g)$ be a compact Riemannian manifold, then for any function $\phi:M\rightarrow\mathbb{R}$, and for any small enough  $s>0$ the following is also a Riemannian metric on $M$ 
    \begin{equation}\label{eq:201}
        g_s=(1+s\phi)g\,,
    \end{equation}
    called the first-order conformal deformation of $g$.
\end{definition}
The variation of sectional curvature of a metric under the first order conformal deformation is given by the following lemma \cite{[Strake]}; cf. \cite[Chapter 3, Corollary 3.4]{[Be2]}.

\begin{lemma}\label{lem:1}Let $(M,g)$ be a Riemannian manifold with sectional curvature $\mathrm{sec}_g\geq0$, and let $X,Y\in T_pM$ be $g$-orthonormal vectors such that $\mathrm{sec}_g(X\wedge Y)=0$. Consider a first-order conformal deformation $g_s=(1+s\phi)g$ of $g$. The first variation of $\mathrm{sec}_{g_s}(X\wedge Y)$ is\begin{equation}\label{eq:202}
        \frac{\mathrm{d}}{\mathrm{d}s}\mathrm{sec}_{g_s}(X\wedge Y)\vert_{s=0}=-\frac{1}{2}\mathrm{Hess}\,\phi(X,X)-\frac{1}{2}\,\mathrm{Hess}\,\phi(Y,Y)\,.
    \end{equation}
\end{lemma}

We will also need the following elementary fact 
\cite[Chapter 3, Lemma 3.5]{[Be2]}.
\begin{lemma}\label{lem:2}
    Let $f:[0,S]\times K\rightarrow \mathbb{R}$ be a smooth function, where $S>0$ and $K$ is a compact subset of a manifold. Assume that $f(0,x)\geq0$ for all $x\in K$, and $\frac{\partial f}{\partial s}>0$ if $f(0,x)=0$. Then there exists $s_*>0$ such that $f(s,x)>0$ for all $x\in K$ and $0<s<s_*$.
\end{lemma}

Wilking's metric $(S^2\times S^3, g_W)$ has positive sectional curvature away from a hypersurface $Z$; see discussion at the end of Section \ref{Section Wilking's Metric}. The biorthogonal and distance curvatures are positive inside $Z$ except for points that lie in four disjoint copies of $S^2$. Every point in these four 2-spheres carries an $S^1$ worth of flat 2-planes. Denote these four 2-spheres by\begin{equation}\label{2-spheres with flats}\{S^2_i : i = 1, 2, 3, 4\}\end{equation}We only deform Wilking's metric near these four submanifolds. Let\begin{equation}\chi_i: S^2\times S^3 \rightarrow \mathbb{R}\end{equation}denote a bump function of $S^2_i$, i.e., a nonnegative function that is identically zero outside a tubular neighborhood of $S^2_i$, and identically one in a smaller tubular neighborhood of $S^2_i$. Finally, we define four functions\begin{equation}\{\psi_i:S^2\times S^3 \rightarrow \mathbb{R}: i = 1, 2, 3, 4\}\end{equation} as\begin{equation}\psi_i(p)=\mathrm{dist}_{g_W}(p,S^2_i)^2\end{equation} for $p\in S^2\times S^3$, where $\mathrm{dist}_{g_W}$ is the metric distance function on $(S^2\times S^3, g_W)$. Let $\phi:S^2\times S^3 \rightarrow \mathbb{R}$ be a function defined as
    \begin{equation}\label{eq:203}
        \phi=-\chi_1\psi_1-\chi_2\psi_2-\chi_3\psi_3-\chi_4\psi_4\,,
    \end{equation}
    and consider the first-order conformal deformation of $g_W$ given by
    \begin{equation}\label{eq:204}
        g_s=(1+s\phi)g_W\,.
    \end{equation}
    Note that at a point $p\in S^2_i$ we have 
    \begin{equation}\label{eq:205}
        \mathrm{Hess}\,\phi(X, X)=-\mathrm{Hess}\,\psi_i(X, X)=-2g_W(X_\perp, X_\perp)^2=-2\|X_\perp\|^2_{g_W} \,,
    \end{equation}
    where $X_\perp$ denotes the component of $X$ perpendicular to $S^2_i$.
For each $\theta>0$, consider the compact subset of $(S^2\times S^3)\times \mathrm{Gr_2}(T(S^2\times S^3))\times \mathrm{Gr_2}(T(S^2\times S^3))$ given by\begin{equation}K_\theta:=\{(p,\sigma,\sigma'):\sigma,\sigma'\in\mathrm{Gr}_2(T_p(S^2\times S^3)),\mathrm{dist}(\sigma,\sigma')\geq\theta\},\end{equation} and define
    \begin{equation}\label{eq:206}
        \begin{split}
            f&:[0,S]\times K_\theta\rightarrow \mathbb{R}\, \\
            f(s,(p,\sigma,\sigma'))&:=\frac{1}{2}(\mathrm{sec}_{g_s}(\sigma)+\mathrm{sec}_{g_s}(\sigma'))\,.
        \end{split}
    \end{equation}
Notice that $f(0,(p,\sigma,\sigma'))\geq 0$ since $\mathrm{sec}_{g_s}\geq 0$. Furthermore, $f(0,(p,\sigma,\sigma'))=0$ only for\begin{equation}p\in S^2_1\cup S^2_2 \cup S^2_3\cup S^2_4\end{equation}since these are the only points of $S^2\times S^3$ that have vanishing biorthogonal and distance curvatures. Let $(p,\sigma,\sigma')$ be such that $f(0,(p,\sigma,\sigma'))=0$ and let $\sigma=X\wedge Y$ and $\sigma'=Z\wedge W$, where $X, Y$ are $g_W$-orthonormal, and $Z, W$ are $g_W$-orthonormal. Then, by Lemma \ref{lem:1} and equation \eqref{eq:205}, at these points of $K_\theta$ we have
    \begin{equation}\label{eq:68}
        \begin{split}
            \frac{\partial f}{\partial s}\vert_{s=0}&=\frac{\mathrm{d}}{\mathrm{d}s}(\mathrm{sec}_{g_s}(X\wedge Y)+\mathrm{sec}_{g_s}(Z\wedge W))\vert_{s=0}= \\
            &=-\frac{1}{2}\mathrm{Hess}\,\phi(X,X)-\frac{1}{2}\mathrm{Hess}\,\phi(Y,Y)-\frac{1}{2}\mathrm{Hess}\,\phi(Z,Z)-\frac{1}{2}\mathrm{Hess}\,\phi(W,W)= \\
            &=\|X_\perp\|^2_{g_W}+\|Y_\perp\|^2_{g_W}+\|Z_\perp\|^2_{g_W}+\|W_\perp\|^2_{g_W}>0\,.
        \end{split}
    \end{equation}
    The previous expression is strictly greater than zero. Indeed, since $X\wedge Y$ and $Z\wedge W$ are different 2-planes, $\mathrm{span}\{X,Y,Z,W\}$ is at least three-dimensional while the submanifolds (\ref{2-spheres with flats}) are two-dimensional. Hence, at least one of the perpendicular components $X_\perp,Y_\perp,Z_\perp,W_\perp$ is nonzero and (\ref{eq:68}) is greater than zero. Since the assumptions of Lemma \ref{lem:2} for the function (\ref{eq:206}) are satisfied, we conclude that there is an $s_*$ such that $f(s,(p,\sigma,\sigma'))>0$ for all $(p,\sigma,\sigma')\in K_\theta$ and $0<s<s_*$. This is precisely the condition $\mathrm{sec}_{g_s}^\theta>0$ of Item (a) of Theorem \ref{Theorem A}. The claims of Item (b) and Item (c) follow from our construction; cf. \cite{[Be1]}. The claim of Item (d) follows from \cite[Proposition 4.1]{[Be1]}.  As Bettiol observed in his construction of metrics of positive distance curvature on $S^2\times S^2$ \cite[Section 4.4]{[Be1]}, for every $\theta > 0$, there are 2-planes in $(S^2\times S^3, g^{\theta})$ with negative sectional curvature. This completes the proof of Theorem \ref{Theorem A}.\hfill $\square$

\begin{remark}The metrics $(S^2\times S^3, g^\theta)$ of positive distance curvature can be made invariant under the action of certain Deck transformations including the product $\Z/2\oplus \Z/2$-action. Indeed, it is possible to perform a local conformal deformation on the orbit space $(\R P^2\times \R P^3, g_W)$ equipped with Wilking's metric of almost positive curvature, and a similar statement to Theorem \ref{Theorem A} holds for $(\R P^2\times \R P^3, g^\theta)$; cf. \cite[Section 4.6]{[Be1]}.

\end{remark}

\subsection{Proof of Corollary \ref{Corollary B}} We will use a case of the classification up to diffeomorphism of simply connected 5-manifolds with vanishing second Stiefel-Whitney class due to Smale \cite[Theorem A]{[S]}.

\begin{theorem} A closed simply connected 5-manifold $M$ with torsion-free homology $H_2(M; \Z) = \Z^k$ and zero second Stiefel-Whitney class $w_2(M) = 0$ is determined up to diffeomorphism by its second Betti number $b_2(M)$. In particular, $M$ is diffeomorphic to a connected sum\begin{equation}\label{5-mflds}\{S^5\#k(S^2\times S^3): k = b_2(M)\}.\end{equation}\end{theorem} Theorem \ref{Theorem A} and Bettiol's result regarding the positivity of biorthogonal curvature under connected sums \cite[Proposition 7.11]{[Be2]} imply that every 5-manifold in the set (\ref{5-mflds}) admits a Riemannian metric of positive biorthogonal curvature. 
\hfill $\square$

\begin{remark}\label{Remark 2}It is natural to ask if the hypothesis $w_2(M) = 0$ of Corollary \ref{Corollary B} can be removed. Barden has shown that a closed simply connected 5-manifold with torsion-free second homology group is diffeomorphic to a connected sum of copies of $S^2\times S^3$ and the total space $S^3\widetilde{\times} S^2$ of the nontrivial $3$-sphere bundle over the $2$-sphere \cite{[Ba]}. It is currently unknown if there is a metric of almost positive sectional curvature on $S^3\widetilde{\times} S^2$. Unlike $S^2\times S^3$, the nontrivial bundle does not arise as a biquotient that satisfies the symmetry hypothesis needed to apply Wilking's doubling trick; see DeVito's classification of free circle actions on $S^3\times S^3$ in \cite{[DeV]}; cf. \cite{[DeVito1], [DeVito2]}.
\end{remark}

\subsection{Proof of Proposition \ref{Proposition Wu manifold}}The symmetric space metric on $\SU(3)/\SO(3)$ is the metric that makes the canonical surjection 
\begin{equation}\label{eq:1}
    \begin{split}
        \pi :\SU(3) &\to \SU(3)/\SO(3) \\
                          u &\mapsto u\SO(3)\,,
    \end{split}
\end{equation}
into a Riemannian submersion, where $\SU(3)$ is equipped with a bi-invariant metric. The left action of $\SU(3)$ on $\SU(3)/\SO(3)$ induced from the left multiplication on $\SU(3)$ by \eqref{eq:1} is transitive and isometric for the symmetric space metric. This means that we can study curvature at one point of $\SU(3)/\SO(3)$ and isometrically translate the results to any other point. The Cartan decomposition that corresponds to $\SU(3)/\SO(3)$ 
\begin{equation}\label{eq:2}
    T_e\SU(3) \simeq \mathfrak{su}(3) = \mathfrak{so}(3) \oplus \mathfrak{so}(3)^\perp \,.
\end{equation} is orthogonal with respect to the bi-invariant metric and it is precisely the decomposition of $T_e\SU(3)$ into vertical and horizontal subspaces of the Riemannian submersion \eqref{eq:1}. Hence, we have 
\begin{equation}\label{eq:3}
    T_{ \SO(3) }( \SU(3)/\SO(3) ) \simeq \mathfrak{so}(3)^\perp\,.
\end{equation}
To conclude that $ \mathrm{SU}(3)/\mathrm{SO}(3)$ has positive biorthogonal curvature, we need to show that no two flat 2-planes are orthogonal to each other. A result of Tapp \cite[Theorem 1.1]{[Tapp]} implies that a 2-plane on $ \mathrm{SU}(3)/\mathrm{SO}(3) $ is flat if and only if its horizontal lift is flat. Thus, it is enough to consider horizontal flat 2-planes at the identity of $ \mathrm{SU}(3) $. 

A horizontal 2-plane $ X \wedge Y \subset \mathfrak{so}(3)^\perp $ at the identity of $\mathrm{SU}(3)$ is flat if and only if $ [ X,Y ] = 0 $. Since the maximal number of linearly independent commuting matrices in $\mathfrak{su}(3)$ is two, every horizontal flat 2-plane corresponds to a maximal abelian subalgebra of $\mathfrak{so}(3)^\perp $
\begin{equation}\label{eq:4}
    \mathrm{span}_\R \{ X, Y \} = \mathfrak{a}_0 \subset \mathfrak{so}(3)^\perp\,. 
\end{equation}
By a fundamental fact about Cartan decomposition, see \cite[Proposition 7.29]{[Knapp]} for the precise statement, any two maximal abelian subalgebras of $\mathfrak{so}(3)^\perp$ are conjugate by an element of $\SO(3)$. This means that by fixing one maximal abelian subalgebra, or one horizontal flat 2-plane we can parametrize all horizontal flat 2-planes by $\SO(3)$. In what follows we will obtain an explicit parametrization of horizontal flat 2-planes at the identity of $\SU(3)$, and so a parametrization of flat 2-planes at a point of $\SU(3)/\SO(3)$ by choosing a basis for $ \mathfrak{su}(3) $, fixing a horizontal flat 2-plane and parametrizing $ \SO(3) $ by Euler angles. We use this explicit parametrization to show that no two flat 2-planes can be orthogonal.
For the basis of $ \mathfrak{su}(3)$, we choose $ \{-i\lambda_i\}_{i=1,...,8} $, where the $ \lambda_i $'s are traceless, self-adjoint 3 by 3 matrices known as the Gell-Mann matrices \cite{[Gell]}. The scalar product on $\mathfrak{su}(3)$ that corresponds to the bi-invariant metric is 
\begin{equation}\label{eq:5}
   \langle X, Y \rangle = -\frac{1}{2}\mathrm{Tr}(XY)\,, 
\end{equation}
for  $X, Y \in \mathfrak{su}(3)$ and the basis $ \{-i\lambda_i\}_{i=1,...,8} $ is orthonormal with respect to (\ref{eq:5}). The Cartan decomposition \eqref{eq:2} in this basis is 
\begin{equation}\label{eq:6}
    \mathfrak{so}(3) = \mathrm{span}_\R\{ -i\lambda_2, -i\lambda_5, -i\lambda_7\}
\end{equation}and
\begin{equation}\label{eq:7}
    \mathfrak{so}(3)^\perp = \mathrm{span}_\R\{ -i\lambda_1, -i\lambda_3, -i\lambda_4,-i\lambda_6,-i\lambda_8\}\,.
\end{equation}
 Matrices $\lambda_3$ and $\lambda_8$ are diagonal, so we use $-\lambda_3 \wedge \lambda_8$ for the reference horizontal flat 2-plane. Every horizontal flat 2-plane, $ X \wedge Y $, with $ X, Y \in \mathfrak{so}(3)^\perp $ such that $ [X,Y]=0 $, can now be written as
\begin{equation}\label{eq:8}
    X \wedge Y = -\mathrm{Ad}_r (\lambda_3 \wedge \lambda_8 )\,,
\end{equation}
for some $ r \in \mathrm{SO}(3) $. Suppose that $ X\wedge Y $ and $ X'\wedge Y' $ are two such 2-planes with, $X\wedge Y$ given by \eqref{eq:8}, and $ X' \wedge Y' $ by
\begin{equation}\label{eq:9}
    X' \wedge Y' = -\mathrm{Ad}_{r'} (\lambda_3 \wedge \lambda_8 )\,,
\end{equation}
for some $ r' \in \mathrm{SO}(3) $. For the 2-planes \eqref{eq:8} and \eqref{eq:9} to be orthogonal it is necessary and sufficient that the equations
\begin{equation}\label{eq:10}
    \langle \mathrm{Ad}_r \lambda_3, \mathrm{Ad}_{r'} \lambda_3 \rangle = 0\,,
\end{equation}
\begin{equation}\label{eq:11}
    \langle \mathrm{Ad}_r \lambda_3, \mathrm{Ad}_{r'} \lambda_8 \rangle = 0\,,
\end{equation}
\begin{equation}\label{eq:12}
    \langle \mathrm{Ad}_r \lambda_8, \mathrm{Ad}_{r'} \lambda_3 \rangle = 0\,,
\end{equation}and
\begin{equation}\label{eq:13}
    \langle \mathrm{Ad}_r \lambda_8, \mathrm{Ad}_{r'} \lambda_8 \rangle = 0\
\end{equation}hold.
Using the $ \mathrm{Ad} $-invariance of the bi-invariant metric, equations \eqref{eq:10}, \eqref{eq:11}, \eqref{eq:12}, and \eqref{eq:13} can be rewritten as
\begin{equation}\label{eq:25}
    \langle \lambda_3, \mathrm{Ad}_{r^{-1}r'} \lambda_3 \rangle = 0\,,
\end{equation}
\begin{equation}\label{eq:26}
    \langle \lambda_3, \mathrm{Ad}_{r^{-1}r'} \lambda_8 \rangle = 0\,,
\end{equation}
\begin{equation}\label{eq:27}
    \langle \lambda_8, \mathrm{Ad}_{r^{-1}r'} \lambda_3 \rangle = 0\,,
\end{equation}and
\begin{equation}\label{eq:28}
    \langle \lambda_8, \mathrm{Ad}_{r^{-1}r'} \lambda_8 \rangle = 0\,.
\end{equation}
We now use the Euler angle parametrization of $ \mathrm{SO}(3) $ to write $ r^{-1}r' \in \mathrm{SO}(3) $ as
\begin{equation}\label{eq:29}
    r^{-1}r' = \mathrm{exp}(- i \lambda_2 x) \mathrm{exp}(- i \lambda_5 y) \mathrm{exp}(- i \lambda_2 z)\,,
\end{equation}
where $ x,y,z \in \mathbb{R}  $. Plugging \eqref{eq:29} into equations \eqref{eq:25}, \eqref{eq:26}, \eqref{eq:27}, and \eqref{eq:28} and calculating the traces explicitly, we find 
\begin{equation}\label{eq:30}
    0 = \langle  \lambda_3, \mathrm{Ad}_{r^{-1} r'} \lambda_3 \rangle = \frac{1}{4} \mathrm{cos}(2x) \left (3 + \mathrm{cos}(2y) \right ) \mathrm{cos}(2z) - \mathrm{sin}(2x)  \mathrm{cos}(y) \mathrm{sin}(2z)\,,
\end{equation}
\begin{equation}\label{eq:31}
    0 = \langle  \lambda_3, \mathrm{Ad}_{r^{-1} r'} \lambda_8 \rangle = - \frac{\sqrt{3}}{2} \mathrm{cos}(2x)\mathrm{sin}^2(y)\,,
\end{equation}
\begin{equation}\label{eq:32}
     0 = \langle  \lambda_8, \mathrm{Ad}_{r^{-1} r'} \lambda_3 \rangle = - \frac{\sqrt{3}}{2} \mathrm{cos}(2z)\mathrm{sin}^2(y)\,,
\end{equation}and
\begin{equation}\label{eq:33}
     0 = \langle  \lambda_8, \mathrm{Ad}_{r^{-1} r'} \lambda_8 \rangle =\frac{1}{4}(1+3\mathrm{cos}(2y))\,.
\end{equation}
Equations \eqref{eq:31}, \eqref{eq:32}, and \eqref{eq:33} imply $\mathrm{cos}^2(y)=1/3$ and $\mathrm{cos}(2x)=\mathrm{cos}(2z)=0$. Plugging this into equation \eqref{eq:30}, we obtain\begin{equation} \langle  \lambda_3, \mathrm{Ad}_{r^{-1} r'} \lambda_3 \rangle \neq 0,\end{equation} and conclude that there is no solution to the system given by equations \eqref{eq:30}, \eqref{eq:31}, \eqref{eq:32}, and \eqref{eq:33}. This shows that no two 2-flat planes are orthogonal. 
\hfill $\square$\\

\end{document}